\documentclass[11pt]{article}

\usepackage{amssymb,amsthm,amsmath}
\usepackage{algorithm}
\usepackage{algpseudocode}
\usepackage{mathtools}

\usepackage{subfig}
\usepackage[usenames,dvipsnames]{xcolor}
\usepackage[colorlinks,citecolor=blue,linkcolor=BrickRed]{hyperref}\usepackage[colorlinks,citecolor=blue,linkcolor=BrickRed]{hyperref}
\usepackage{fullpage}
\usepackage{graphicx}
\usepackage{xspace}
\usepackage{enumitem}
\usepackage{thmtools,thm-restate}
\usepackage{wrapfig}
\usepackage{comment}
\usepackage{cleveref}

\graphicspath{{figures/}}

\newcommand{\ip}[2]{\langle #1, #2 \rangle}
\newcommand{\supp}{\mathsf{supp}}

\newtheorem{theorem}{Theorem}[section]

\newtheorem{lemma}[theorem]{Lemma}

\newtheorem{remark}[theorem]{Remark}
\newtheorem{claim}[theorem]{Claim}

\newif\ifFULL
\FULLtrue

\newcommand{\IGNORE}[1]{}

\usepackage{tikz}
\usetikzlibrary{arrows}
\usetikzlibrary{arrows.meta}
\usetikzlibrary{shapes}
\usetikzlibrary{backgrounds}
\usetikzlibrary{positioning}
\usetikzlibrary{decorations.markings}
\usetikzlibrary{patterns}
\usetikzlibrary{calc}
\usetikzlibrary{fit}
\usetikzlibrary{snakes}
\tikzset{
    >=stealth',
    pil/.style={
           ->,
           thick,
           shorten <=2pt,
           shorten >=2pt,}
}
\tikzset{->-/.style={decoration={
  markings,
  mark=at position .5 with {\arrow{>}}},postaction={decorate}}}

\newcommand{\I}{\mathbb{I}}

\newcommand{\p}{\ensuremath{\mathbb{P}}}

\newcommand{\E}{\mathbb{E}}

\newcommand{\ind}{\ensuremath{\mathbf{1}}}
\newcommand{\eat}[1]{}
\newcommand{\hide}[1]{{\Large \color{red} Contents here are hidden! To reveal contents, remove this command.}}

\newcommand{\disc}{\ensuremath{\mathrm{disc}}}

\newcommand{\poly}{\ensuremath{\mathrm{poly}}}

\newcommand{\CE}{\mathcal{E}}
\newcommand{\CG}{\mathcal{G}}

\newcommand{\CU}{\mathcal U}
\newcommand{\CN}{\mathcal N}
\renewcommand{\Pr}{\p}
\newcommand{\pmone}{\{\pm 1\}}

\allowdisplaybreaks

{\hspace*{\fill}$\Box$\par}

\title{Smoothed Analysis of the Koml\'{o}s Conjecture}

\author{ Nikhil Bansal\thanks{University of Michigan, Ann Arbor. \texttt{bansal@gmail.com}.} \and
 Haotian Jiang\thanks{University of Washington, Seattle. \texttt{jhtdavid@cs.washington.edu}.
 } \and
  Raghu Meka\thanks{University of California, Los Angeles. \texttt{raghum@cs.ucla.edu}. 
 } \and
 Sahil Singla\thanks{Georgia Institute of Technology, Atlanta. \texttt{ssingla@gatech.edu}. } \and 
  Makrand Sinha\thanks{Simons Institute and University of California, Berkeley. \texttt{makrand@berkeley.edu}.}
 }

\date{}

\begin{document}

  \maketitle
  
\begin{abstract}

The well-known Koml\'os conjecture states that given $n$ vectors in $\mathbb{R}^d$ with Euclidean norm at most one, there always exists a $\pm 1$ coloring such that the $\ell_{\infty}$ norm of the signed-sum vector is a constant independent of $n$ and $d$. We prove this conjecture in a smoothed analysis setting where the vectors are perturbed by adding a small Gaussian noise and when the number of vectors $n =\omega(d\log d)$. The dependence of $n$ on $d$ is the best possible even in a completely random setting. 

Our proof relies on a weighted second moment method, where instead of considering uniformly randomly colorings we apply the second moment method on an implicit distribution on colorings obtained by applying the Gram-Schmidt walk algorithm to a suitable set of vectors. The main technical idea is to use various properties of these colorings, including subgaussianity, to control the second moment.

  \end{abstract}


\section{Introduction}
A central question in discrepancy theory is the following Koml\'{o}s problem: given vectors $v_1,\ldots,v_n \in \mathbb{R}^d$ with Euclidean length at most $1$, i.e.,~$\|v_i\|_2 \leq 1$ for all $i \in [n]$, find signs $x_i \in \{-1,1\}$ for $i\in [n]$ to minimize the discrepancy $\| \sum_{i=1}^n x_i v_i \|_\infty$. 
The long-standing Koml\'{o}s conjecture says that 
the discrepancy of any collection of such vectors is $O(1)$, independent of $n$ and $d$.
An important special case (up to scaling by $t^{1/2}$) is the Beck-Fiala problem, where the vectors $v_1,\ldots,v_n \in \{0,1\}^d$  and each $v_i$ has at most $t$ ones, so $\|v_i\|_2 \leq t^{1/2}$. 
Here, the Koml\'{o}s conjecture reduces to the Beck-Fiala conjecture \cite{BeckFiala-DAM81}, which says that the discrepancy is $O(t^{1/2})$. 
The question of either proving or disproving these conjectures has received a lot of attention, and
after a long line of work, the current best bounds for the Komlos and the Beck-Fiala problem are $O((\log n)^{1/2})$ and $O( (t \log n)^{1/2}) $ respectively, due to Banaszczyk \cite{Banaszczyk-Journal98}.

Motivated by the lack of progress for general worst-case instances, there has been a lot of recent work on these problems for random instances, with several interesting results and techniques, see e.g.,~\cite{EzraLovett-Random19,BansalMeka-SODA19,HobergRothvoss-SODA19,FS20, Potukuchi-ICALP20, TMR-COLT20, AW21, macrury2021phase}. 
In this work, we consider the Koml\'{o}s problem in the 
more general setting of 
smoothed analysis, where the input is generated by 
taking an arbitrary worst case Koml\'{o}s instance 
and perturbing it randomly. 
The smoothed analysis model was first introduced by Spielman and Teng~\cite{ST-JACM04}, and it interpolates nicely between worst case and average case analysis, and has been used extensively since then to study various problems. Recently smoothed analysis models have also been considered in discrepancy theory in a few other works \cite{BJMSS-ITCS22, HRS-FOCS21} ---
however the setting and focus of these results is quite different, and in particular they are not directly related to the Koml\'os or Beck-Fiala conjectures.


\medskip
\noindent \textbf{Random instances.}
To put our results in the proper context, we first describe the results on random instances.
In general, these results depend on the different regimes of the parameters $d,n$ and $t$, and we  focus here on the more interesting case of $n \gg d$. 

A natural model for random Beck-Fiala instances is
where each entry is $1$ with probability $p=t/d$, so that each column has $t$ ones in expectation.
In a surprising result,
Hoberg and Rothvoss \cite{HobergRothvoss-SODA19} showed
that $\disc(A) \leq 1$ w.h.p.\footnote{This is much better than the $O(t^{1/2})$ bound in the Beck-Fiala conjecture.}  if
$n = \Omega(d^2 \log d)$.
Independently, Franks and Saks \cite{FS20}
showed that $\disc(A) \leq 2$ w.h.p.
if $n = \Omega(d^3  \log^2 d)$, for a more general class of instances.
Both these results use interesting Fourier analysis based  techniques.

It is not hard to see\footnote{ If we fix any coloring $x$ and consider a random instance, a fixed row has discrepancy $O(1)$ with probability $\approx (pn) ^{-1/2}$, so the probability that each row has discrepancy $O(1)$ is $(pn)^{-\Omega(d)}$. As there are (only) $2^n$ possible colorings, a first moment argument already requires that $2^n (pn)^{-d} = \Omega(1)$.} that $n = \Omega(d \log d)$ is necessary for $O(1)$ discrepancy (provided $p$ is not too small).
An important step towards obtaining this optimal dependence was made by
Potukuchi \cite{P18}, who
showed that $\disc(A) \leq  1$ if $n = \Omega(d \log d)$ for the dense case of $p = 1/2$, using the second moment method.
However, the sparse setting with $p\ll 1$ turns out to be more subtle, and was only recently resolved by Altschuler and Weed \cite{AW21} using a more sophisticated conditional second moment method together with Stein's method of exchangeable pairs. They show that $\disc(A) \leq 1$ w.h.p. for $n \geq  \Omega(d \log d)$, for every $p$.

The case of Gaussian matrices with i.i.d.~$\CN(0,1)$ entries has also been considered, where
Turner, Meka and Rigollet \cite{TMR-COLT20} give almost tight bounds for the entire regime, and 
in particular show that for $n = \Omega(d \log d)$ a discrepancy bound of $1/\poly(d)$ holds. 

\medskip
\noindent \textbf{The smoothed Koml\'os model.} We now define our model formally. The input matrix is of the form $A = M+R$, where $M \in \mathbb{R}^{d \times n}$ is some worst-case matrix with columns of $\ell_2$-norm at most $1$ and $R \in \mathbb{R}^{d \times n}$ is a random matrix 
with i.i.d. Gaussian entries distributed as $\mathcal{N}(0,\sigma^2/d)$, where $\sigma \leq 1$.  
The $\sigma^2/d$ variance ensures that each column of $R$ has $\ell_2$-norm roughly $\sigma$ (and hence much less than that of $M$).
Our goal is to understand the  discrepancy of $A$.
We will only be interested in showing the existence of a low discrepancy coloring for $A$, and not in algorithmically finding it (this seems far beyond the current techniques).



\subsection{Results and Techniques}

Our main result is the following.
\begin{theorem}[Smoothed Koml\'os]
\label{thm:main}
Let $\sigma>0$ and $n  = \frac{\omega(d \log d)}{\sigma^{4/3}}$. Then with probability $1 - o_d(1)$, the discrepancy of $M + R$ is at most $1/\poly(d)$, where $M \in \mathbb{R}^{d \times n}$ is an arbitrary Koml\'os instance and $R \in \mathbb{R}^{d \times n}$ has i.i.d. $\mathcal{N}(0,\sigma^2/d)$ entries. 
\end{theorem}

An interpretation of \Cref{thm:main} is that any counter-example to the Koml\'os conjecture (if it exists) will be rigid, or have $n \approx d$.
Also, notice that  dependence of $n$ on $d$ in Theorem~\ref{thm:main} is essentially the best possible, as already evident in the very special case of $M = \bf{0}$ and $\sigma=1$, i.e.,~a random matrix with i.i.d. $\CN(0,1)$ entries, where $n = \Omega(d \log d)$ is necessary to achieve $1/\poly(d)$ discrepancy, as discussed earlier.

\begin{remark}
Our proof techniques also give a high probability bound when $n = \Omega_\sigma(d^{1+\epsilon})$ for any constant $\epsilon>0$. However, we do not explore this direction here. It would be interesting to know if the result also holds with high probability when $n=\omega(d \log d)$, as in the (fully) random setting.
\end{remark}

\begin{remark}
The dependence on the noise parameter $\sigma$ in $n=\Omega_d(1/\sigma)$ in Theorem \ref{thm:main} is necessary, otherwise this would imply an $O(1)$ bound for the worst case Koml\'{o}s problem. In particular, 
each row of the random part $R$ must at least have enough $\ell_1$ norm to offset the discrepancy from the worst case part $M$ (which can be $O((\log n)^{1/2})$ given the currently known results).
As each entry of $R$ has magnitude about $\sigma d^{-1/2}$, we thus require $n  = \Omega_d(1/\sigma)$ for each row of $M+R$ to have discrepancy $O(1)$.
 \end{remark}


The proof of \Cref{thm:main} is based on the classical second moment method, however, it requires several additional ideas beyond those used for random instances, to handle the effect of the worst case part and its interplay with the random part.
We describe these briefly next, and discuss them in more detail in Section \ref{sec:overview}.

\smallskip

\begin{itemize}
    \item \textbf{Weighted second moment method.} Instead of applying the second moment method to the uniform distribution on the $2^n$ colorings, we consider a distribution on low-discrepancy colorings for $M$.
This is necessary as for a random coloring $x \in \{-1,1\}^n$, a typical entry of $Mx$ will scale as $\sqrt{n/d}$, which is very unlikely to be {\em  cancelled} by the discrepancy of the random part $Rx$, which typically scales as $\sigma \sqrt{n/d}$ (note that we want to show the existence of some $x$ such that  $(Rx)_i \approx   -(Mx)_i$ for each coordinate $i \in [d]$).

\smallskip

\item \textbf{Subgaussianity of colorings.} To ensure that $\|Mx\|_\infty$ is typically small, we consider the (implicit) distribution on colorings produced by the Gram-Schmidt (GS) algorithm \cite{BansalDGL18} applied to $M$, which ensures that $Mx$ is a $1$-subgaussian vector \cite{hssz19} (details in Section~\ref{sec:overview}). 

However, apriori the GS algorithm only guarantees that $Mx$ is subgaussian, and says nothing about the distribution on the colorings $x$. For instance, it could be that
any two colorings in the support have the first $9n/10$ coordinates identical, and thus look very non-random. This makes the second moment bounds much worse and harder to control.

To handle this, we use a simple but useful trick to ensure that the distribution on the colorings $x$ produced by the GS algorithm is also $O(1)$-subgaussian. Roughly, this allows us to pretend that colorings $x$ in the GS distribution behave randomly.

\smallskip

\item \textbf{Exploiting subgaussianity to get cancellations across rows.} 
Most importantly, due to the worst case part $M$, doing a row by row analysis 
as is typically done in second moment computations for random instances, only works when $n = \Omega(d^2/\sigma^2)$ (details in Section \ref{sec:overview}). Roughly, the problem is that considering each coordinate of $Mx$ separately completely ignores the global properties across the different coordinates that subgaussianity of $Mx$ implies.

To get the optimal dependence of $n$ on $d$, a key conceptual idea is to analyze {\em all} the rows together  and use the subgaussianity of $Mx$ and $x$ carefully to get various cancellations across the different rows in the second moment computation.  Exploiting subgaussianity also leads to various technical difficulties, as subgaussian vectors can differ from fully random Gaussian vectors in various non-trivial ways.
\end{itemize}

\medskip \noindent
\textbf{Notation.} Throughout this paper, $\log$ denotes the natural logarithm. We use the asymptotic notation $\omega(\cdot)$ or $o(\cdot)$ where the growth is always with respect to $d$ --- sometimes to emphasize this dependence we will also write $\omega_d(\cdot)$ or $o_d(\cdot)$. We write $\E_{x \sim \CG}[f(x)]$ to denote the expectation of a function $f$ where $x$ is sampled from the distribution $\CG$ and we abbreviate this to $\E[f(x)]$ when the distribution is clear from the context. For reals $a, b \in \mathbb{R}$, the notation $[a \pm b]$ is used as a shorthand to denote the interval $[a-b, a+b]$. For a set $S \in \mathbb{R}^d$, we write $\delta S = \{ \delta x \mid  x \in S\}$ to denote the $\delta$ scaling of $S$.

\subsection{Overview and Preliminaries}
\label{sec:overview}
We now give a more detailed overview of the proof and the ideas. We also briefly describe the second moment method and some concepts we need such as subgaussianity and properties of the Gram-Schmidt algorithm.

\medskip \noindent \textbf{Second moment method.} The second moment method (e.g. \cite{alonspencer}) is based on the following Paley-Zygmund inequality. For any non-negative random variable $Z$, we have that
\[ \Pr[Z >0] \geq \frac{(\E[Z])^2}{\E[Z^2]}.\]
So, if $\E[Z^2] = (1+o(1)) (\E[Z])^2$, then this implies that $\Pr[Z > 0] \geq 1-o(1)$.

For constraint satisfaction problems, a standard way to use this  to show that most random instances are feasible is by defining $S = S(R)$ as the number of solutions to an instance $R$, and showing that 
\begin{equation}
\label{eq:sec-moment}
    \E_R[S^2] = (1+o(1)) (\E_R[S])^2,
\end{equation}
which gives that  $\Pr_R[S(R) >0] = \Pr_R[S(R) \geq 1] \geq  1-o(1)$.

Let us consider \eqref{eq:sec-moment} more closely and define $S(R,x) = 1$ if $x$ is a valid solution for instance $R$, and $0$ otherwise. 
Then $S(R) = \sum_{x} S(R,x)$ and \eqref{eq:sec-moment} can be written as
\begin{equation}
    \label{eq:sec-mom1}
    \E_R  \left[\Pr_{x,y \sim \CU} \big[\, S(R,x)=1 \, , \, S(R,y)=1\, \big]\right]  = (1+o(1)) \left(\E_R \left[ \Pr_{x \sim \CU} \left[S(R,x)=1 \right] \right] \right)^2,
\end{equation} 
where $\CU$ is the uniform distribution over all inputs.

\medskip \noindent \textbf{Second moment method for smoothed Koml\'{o}s.} Let $\Delta$ denote the desired discrepancy bound. In our setting,
denote by $S(R,x)=1$ that 
$x \in \pmone^n$ is a feasible coloring for the smoothed Koml\'{o}s instance $M+R$, that is, if  $\| (M+R)x \|_\infty \leq \Delta$. Roughly, this condition means that  $Rx = -Mx$ and hence the discrepancy of the random part $R$ cancels that of the worst case part $M$. 

However, if $x$ is chosen uniformly from $\pmone^n$, it is not hard to see that this cannot work. The entries $(Mx)_i$ will be distributed roughly as $\CN(0,m_i^2)$ where $m_i = (\sum_{j} M_{ij}^2)^{1/2}$ is the $\ell_2$-norm of row $i$ of $M$, and in general will be much larger (around $1/\sigma \gg 1$ times) than the entries $(Rx)_i$.

\medskip \noindent \textbf{Weighted second moment.} To allow a reasonable probability of $Rx$ cancelling $Mx$, a natural idea is to consider a distribution that is mostly supported 
over colorings $x$ with low discrepancy on $M$. So, we will show \eqref{eq:sec-mom1} where $x,y$ are sampled from some another suitable distribution $\CG$ instead of the uniform distribution $\CU$.
Similar ideas have also been used in other contexts such as \cite{Achlioptas-peres}.
Notice that this does not affect $Rx$, as for any fixed $x$, the contribution of the random part $(Rx)_i$  is still distributed as $N(0,n\sigma^2/d)$ (over the randomness of $R$).

A natural candidate is the distribution on colorings produced by the Gram-Schmidt (GS) walk algorithm \cite{BansalDGL18}. 
In particular, we use the following result.

\begin{theorem}[\cite{hssz19}]\label{thm:hssz}
Given vectors $v_1,\dots,v_n\in\mathbb{R}^m$ with $\|v_j\|_2 \leq 1$, the Gram-Schmidt walk algorithm outputs a random coloring $x\in\{-1,1\}^n$ such that $\sum_{j=1}^n x_j v_j$ is 1-subgaussian.
\end{theorem}
Recall that a random
vector $Y \in \mathbb{R}^m$ is $\alpha$-subgaussian if for all test vectors $\theta\in \mathbb{R}^{m}$,
\[ \E\big[ \exp(\langle \theta,Y\rangle )  \big] \le  \exp\bigg(\frac{\alpha^2\| \theta \|_2^2}{2}\bigg).\]
Roughly, this means that $Y$ looks like a Gaussian with variance at most $\alpha^2$ in every direction. 

Let $\CG$ denote the (implicit) distribution over the colorings output by the GS walk algorithm. 
For a coloring $x$, let us denote $P_x := \Pr_R[ S(R,x)=1]$ and for two colorings $x$ and $y$, let 
\[P_{x, y} := \Pr_R\big[ S(R,x)=1\,,\,S(R,y)=1 \big].\]
Then changing the order of expectation in \eqref{eq:sec-mom1} and substituting, our goal is to show that 
\begin{equation}
    \label{eq:sec-moment-G}\E_{x,y \sim \CG} [P_{xy}]  = (1+o(1)) \cdot \E_{x,y\sim \CG}  [P_x P_y].
\end{equation}
However, the set of low discrepancy colorings for $M$ and the distribution $\CG$ can be quite complicated and hard to work with. 
Later, we will ensure that $\CG$ is also $O(1)$-subgaussian, which will suffice for our purposes.  Let us first consider \eqref{eq:sec-moment-G} more closely.

\medskip \noindent \textbf{The key computation.}
As $R$ has i.i.d. Gaussian entries, the quantities $P_x$ and $P_{xy}$ can be written in a very clean way. 
In particular, as $(Rx)_i \sim \CN(0,\sigma^2 n/d)$ for any coloring $x$, and the $(Rx)_i$ are independent for $i \in [d]$, we can write 
\[ P_x ~=~  \Pr_R\bigg[ \bigcap_{i=1}^d ((Rx)_i \in  (Mx)_i \pm \Delta ) \bigg] ~=~ \prod_{i=1}^d \Pr_R \big[  g_i \in (Mx)_i \pm \Delta \big],  \] 
where $g_i \sim N(0,\sigma^2 n/d)$ and $g_i$'s are independent. 

Similarly, for any fixed colorings $x$ and $y$, writing $g_i = (Rx)_i$ and $g_i' = (Ry)_i$ we have
\[  P_{xy} = \prod_{i=1}^d  \Pr_R \big[  g_i \in (Mx)_i \pm \Delta \,, \, g'_i \in (My)_i \pm \Delta \big],  \]
where $g_i$ and $g'_i$ are correlated with $\E[g_i g_i'] = \langle x,y \rangle \cdot \sigma^2/d$.

A standard computation of 2-dimensional gaussian probabilities over rectangles (and ignoring some less crucial terms for the discussion here) gives 
\begin{equation}
    \label{eq:goal-separate-row}
\frac{\Pr\big[g_i \in (Mx)_i \pm \Delta\,,\, g_i' \in (My)_i \pm \Delta\big] }{\Pr\big[g_i \in (Mx)_i \pm \Delta\big] \cdot \Pr\big[g_i' \in (My)_i \pm \Delta\big]  }  
\approx \exp\bigg(\frac{d \langle x,y\rangle (Mx)_i (My)_i}{\sigma^2 n^2}\bigg).
\end{equation}
So to prove \eqref{eq:sec-moment-G}, we could try to show that for each $i\in [d]$,
\begin{equation}
\label{eq:key-row-i}
 \E_{x,y\sim \CG}   \bigg[ \frac{d \langle x,y\rangle (Mx)_i (My)_i}{ \sigma^2  n^2} \bigg]= o\bigg(\frac1d\bigg).
\end{equation}

Indeed, as $|\langle x ,y \rangle| \leq n$ and $(Mx)_i, (My)_i$ are typically $O(1)$  (as  $Mx$ and $My$ are subgaussian), setting $n = \omega(d^2/\sigma^2)$ would suffice to complete the second moment proof. However, this does not give us the optimal $d\log d$  dependence.

Next, we sketch the two ideas to obtain the optimal dependence.

\medskip \noindent \textbf{Subgaussianity of the distribution $\CG$.}
If $x$ and $y$ were random colorings, we would typically expect that $|\langle x ,y \rangle| \approx \sqrt{n}$ instead of $n$ above. To achieve this, we apply the GS walk algorithm to
the $(d+n) \times n$ matrix with $M$ in top $d$ rows and $I_n$ in the bottom $n$ rows. (Note that each column still has $O(1)$ length.) This ensures that the resulting distribution $\CG$ on the colorings $x$ is $O(1)$-subgaussian, while ensuring that $Mx$ is also $O(1)$-subgaussian.

\medskip \noindent \textbf{Handling the rows together.}
Next, to exploit the subgaussianity of $Mx$  and $My$, we look at all the rows together in \eqref{eq:key-row-i} and consider
\begin{equation}
    \label{eq:key-sum}
 \sum_i  \E_{x,y\sim \CG}  \bigg[ \frac{d \langle x,y\rangle (Mx)_i (My)_i}{ \sigma^2  n^2} \bigg] =  \E_{x,y\sim \CG} \bigg[ \frac{d \langle x,y\rangle \langle Mx,My \rangle}{ \sigma^2  n^2} \bigg].
\end{equation}

By the subgaussianity of the colorings $x, y$ and discrepancy vectors $Mx,My$, we expect that
$\E_{x,y \sim \CG} |\langle x,y\rangle| \approx \sqrt{n}$  and $\E_{x,y}  |\langle Mx,My \rangle| \approx \sqrt{d}$.
Roughly speaking, this implies that the right side of \eqref{eq:key-sum}
is typically $d^{3/2}/(\sigma^2 n^{3/2})$, and hence $ n \gg d/ \sigma^{4/3}$ suffices. 

The formal argument needs some more care as $\langle x,y\rangle$ and $\langle Mx,My\rangle$ are correlated, and as we need to bound the exponential moment of $d\langle x,y\rangle \langle Mx,My \rangle/(\sigma^2  n^2)$ in  \eqref{eq:goal-separate-row}, instead of the expectation, which gives the additional (necessary) logarithmic factor of $\log d$.


\section{Proof of the Smoothed Koml\'os Conjecture}
\label{sec:proof}








We use a weighted version of the second moment method as mentioned in the proof overview. Let $\mathcal{G}$ be a distribution over coloring that will be specified later. We define the following random variable $S$ which depends only on the randomness of $R$, 
\begin{align*}
    S = S(R) := \E_{x \sim \mathcal{G}}[\mathbf{1}\{\|(M+R) x\|_\infty \leq \Delta \}],  
\end{align*}
for some parameter $\Delta = 1/\poly(d)$ to be chosen later. 
The purpose of this variable is that the event $\{S > 0\}$ implies there exists a coloring $x \in \mathsf{supp}(\mathcal{G})$ with discrepancy at most $\Delta$. Our goal is to show that $\p(S > 0) = 1 - o(1)$. As explained in the proof overview, this would follow from the Paley-Zygmund inequality if we can establish that
the first moment $\E_{R}[S]$ is always positive, and the second moment satisfies  $\E_R[S^2] = (1 + o(1)) \cdot (\E_R[S])^2$. We next compute the moments.



\medskip
\noindent \textbf{First moment computation.}
We can compute
\begin{align*}
\E_R[S] ~=~ \E_{x \sim \mathcal{G}} \E_{R}[\mathbf{1}\{\|(M+R) x\|_\infty \leq \Delta \}] ~>~ 0,
\end{align*}
where the  strict inequality follows because fixing any outcome $x \sim \mathcal{G}$, the event $\{\|(M+R) x\|_\infty \leq \Delta \}$ happens with positive probability (recall that $R$ is a Gaussian random matrix with each entry $\mathcal{N}(0,\sigma^2/d)$).

\medskip 
\noindent \textbf{Second moment computation.}
For any $i \in [d]$, denote by $m_i$ and $r_i$ the $i^\text{th}$ row of the matrices $M$ and $R$ respectively. The second moment is given by 
\begin{align*}
\E_R[S^2] 
& = \E_R\left[\E_x[\mathbf{1}\{\|(M+R) x\|_\infty \leq \Delta \}] \cdot \E_y[\mathbf{1}\{\|(M+R) y\|_\infty \leq \Delta \}] \right] \\
& = \E_R \E_{x,y}\left[\mathbf{1}\left\{\|(M+R) x\|_\infty \leq \Delta, \|(M+R) y\|_\infty \leq \Delta \right \} \right] \\
& = \E_{x,y} \left[ \p_R \left(\|(M+R) x\|_\infty \leq \Delta, \|(M+R) y\|_\infty \leq \Delta  \right) \right] =  \E_{xy}[P_{xy}],
\end{align*}
where we define \[P_{x,y} := \p_R(\|(M+R)x\|_\infty\leq \Delta,\|(M+R)y\|_\infty\leq \Delta ).\] Similarly, denoting \[P_x := \p_R(\|(M+R)x\|_\infty\leq \Delta),\] we also have 
\begin{align*}
(\E_R[S])^2 &= (\E_x \p_R(\|(M+R) x\|_\infty \leq \Delta)) \cdot (\E_y \p_R(\|(M+R) y\|_\infty \leq \Delta)) \\
& = \E_{x,y} \left[\p_R(\|(M+R) x\|_\infty \leq \Delta) \cdot \p_R(\|(M+R) y\|_\infty \leq \Delta) \right] \\
\ & = \E_{xy}[P_x \cdot P_y].
\end{align*}

To compare the quantities $\E_R[S^2]$ and $(\E_R[S])^2$, we first consider a distribution over colorings. A natural distribution to consider is the distribution on colorings derived from the Gram-Schmidt walk which ensures that the discrepancy vector $Mx$ is $1$-subgaussian if $x$ is sampled from this distribution. However, we shall also need that the colorings  $x$ themselves have a subgaussian tail as well as some additional nice properties that will be useful to compute the second moment. In particular, we prove the following lemma in \Cref{sec:Truncated}. 


\begin{lemma}[Truncated Gram-Schmidt Distribution]\label{lem:gram-schmidt}
Let $M \in \mathbb{R}^{ d \times n}$ be a worst-case Koml\'os instance. Then, for any constant $C' > 1$ there exists a distribution $\mathcal{G}$ over colorings $x \in \{\pm 1\}^n$ satisfying the following properties:
\begin{itemize}
    \item Almost Constant Euclidean Norm for the discrepancy vectors: for every $x \in \supp(\mathcal{G})$ , we have $\|Mx\|_2 \in [r \pm \Delta]$ where $r = O(d^{1/2})$ and $\Delta = d^{-C'}$.
    \item Almost subgaussian tails for the colorings and discrepancy vectors: there exists a constant $C$ depending on $C'$, such that for every $u \in \mathbb{S}^{n-1}$, 
        \[ \Pr_{x \sim \CG} \left[|\ip{x}{u}| \ge t\right]  \le 2d^C \cdot e^{-t^2/8} \text{ and } 
        \Pr_{x\sim \CG}\left[|\ip{M x}{u}| \ge t\right] 
        \le 2d^C\cdot e^{-t^2/8}. \]
\end{itemize}
\end{lemma}

Since the colorings sampled from the above distribution are subgaussian, $|\ip{x}{y}| \le n/2$ holds with high probability. To compute the second moment to a good precision, we need a careful comparison of the ratio ${P_{xy}}/({P_x \cdot P_y})$ for any two colorings $x$ and $y$ where this event occurs. We show the following bound in this case (proof in \Cref{sec:claims}).


\begin{claim}[Strong bound] \label{claim:second_moment_small_epsilon}
For any two colorings $x,y \in \supp(\CG)$, denote $\epsilon = \epsilon(x, y) = {\ip{x}{y}}/n$. If $|\epsilon| \le 1/2$, then we have 
\begin{align*} 
{P_{x,y}} \leq  P_x P_y \cdot \beta(x,y) \text{ where } \beta(x,y) = \exp\left(\delta_1 + d \epsilon^2 + {d \delta^2 \epsilon^2} + {\delta^2 \epsilon} \cdot \langle Mx, My \rangle \right),
\end{align*}
where the scaling factor $\delta := \frac{\sqrt{d}}{\sigma \sqrt{n}}$ and the error parameter $\delta_1 \le 1/\poly(d)$.  
\end{claim}


When the low probability event $|\epsilon|\ge 1/2$ occurs, we use the weak bound $P_{xy} \le \min\{P_x, P_y\}$.

As $x$ is sampled from the truncated Gram-Schmidt distribution, the probabilities $P_x$ turn out to be almost constant for all colorings  $x \in \supp(\mathcal{G})$ as the following claim  shows.

\begin{claim} \label{claim:large_epsilon}
For any coloring $x \in \supp(\CG)$, 
\begin{align} 
 P_x = \exp(\delta_x) \cdot  p \text{ where } p := \left(\frac{\delta \Delta}{\sqrt{2 \pi}}\right)^d \exp\left(-\frac{\delta^2 r^2}{2} \right),
\end{align}
with the scaling factor $\delta := \frac{\sqrt{d}}{\sigma \sqrt{n}}$ and the error parameter $\delta_x$ satisfying $|\delta_x| \le \delta_1 \le 1/\poly(d)$.  
\end{claim}
The proof of this claim is in \Cref{sec:claim}.

We now focus on the case when $|\epsilon| \leq 1/2$.  
When we take $x,y \sim \CG$, as $P_x$ and $P_y$ are essentially constant, by \Cref{claim:second_moment_small_epsilon}, applying the second moment method reduces  to bounding $\beta(x,y)$ as defined in \Cref{claim:second_moment_small_epsilon}. To do this, we will use the properties,
as described in \Cref{lem:gram-schmidt}, 
of the underlying random variables $x$ and $Mx$. 
The following technical lemma gives a  bound on the exponential moment for such random variables. 


\begin{lemma} \label{lem:truncated_tail}
Let $X$ be a non-negative random variable $X$ that satisfies 
\begin{align*}
    \p(X \geq t) \leq d^{C_1} \cdot e^{ - t^2/8} \text{ for any } t > 0,
\end{align*}
for some fixed constant $C_1 > 0$.
Then for any $\lambda = c_2\sqrt{\log d}$ with $c_2 \ge \sqrt{32C_1}$, 
\begin{align*}
    \E[\exp(X^2/\lambda^2)] \leq 1 + 32C_1/c_2^2 + o_d(1).
\end{align*}
\end{lemma}

We shall prove this lemma in \Cref{sec:Truncated}. 

We can now complete the proof of \Cref{thm:main} by comparing $\E_R[S^2]$ and $(\E_R[S])^2$. We show that
\begin{lemma}\label{lem:second-moment}
For $n =\omega(d\log d) \sigma^{-4/3}$, we have 
    $$(\E_R[S])^2 = p^2 (1 - o_d(1)) \text{ and }  \E_R[S^2] = p^2 (1 + o_d(1)).$$
\end{lemma}

The above implies that $\E_R[S^2] = (1+o(1)) (\E_R[S])^2$, and thus the Paley-Zygmund inequality implies \Cref{thm:main} as discussed in the proof overview. 

\begin{proof}[Proof of \Cref{lem:second-moment}]
    For the first moment, \Cref{claim:large_epsilon} implies that
    \[ (\E_R[S])^2 = \E_{x,y \sim \CG}[P_x P_y] = p^2\E[\exp(\delta_x + \delta_y)] \ge p^2 \exp(-2\delta_1).\]
    Since $0< \delta_1 \le 1/\poly(d)$, the bound follows.
    
    To compute the second moment, $\E_R[S^2]  = \E_{x,y \sim \CG}[P_{xy}]$, we define $\CE$ to be the event that the colorings $x, y \sim \CG$ satisfy $|\ip{x}{y}| > n/2$ and compute the contribution to the expectation under $\CE$ and its complement separately. In particular, 
    using \Cref{claim:large_epsilon} and \Cref{claim:second_moment_small_epsilon}, we have
\begin{align}\label{eqn:main}
 \E_R[S^2] &= \E_{x,y \sim \CG}[P_{x,y}] \leq \p_{x,y \sim \CG}[\CE] \cdot p +  \E_{x,y \sim \CG}\big[P_x P_y \beta(x,y) \cdot \ind[\overline\CE] \big]
\end{align}
For the first term in \eqref{eqn:main}, since $n \ge d$, \Cref{lem:gram-schmidt} implies that 
\[ \p_{x, y \sim \CG}[\CE] \le \poly(d) \cdot e^{-n/4} \le e^{-n/8}.\]
Thus, using the exact bound for $p$ from \Cref{claim:large_epsilon} and that $\delta \Delta \le \poly(\sigma/d)$, the first term 
\begin{align}\label{eqn:first}
    \p_{x,y \sim \CG}[\CE] \cdot p =  p^2 \cdot \p[\CE] \cdot p^{-1} &\le p^2 \cdot e^{-n/8} \cdot \left(\frac{\sqrt{2 \pi}}{\delta \Delta} \right)^d \exp\left(\frac{\delta^2 r}{2}\right) \notag \\
    \ & \leq p^2 \cdot e^{-n/8} \cdot \exp\left(O\left(d \log (dn/\sigma) + d^2/(\sigma^2 n) \right)\right) \notag \\
    & = p^2 \cdot o_d(1),
\end{align}
when $n =\omega(d\log d)\sigma^{-4/3}$. In particular, as $\sigma \leq 1$, we have
$n/8  \gg d\log (dn/\sigma) + d^2/(\sigma^2n)$.

For the second term in \eqref{eqn:main}, using \Cref{claim:large_epsilon}, we have that $P_x = p \cdot \exp(\delta_x)$ where $|\delta_x| \le | \delta_1| \le 1/\poly(d)$. Thus, 
\begin{align}\label{eqn:1}
    \E_{x,y \sim \CG}\big[P_x P_y \beta(x,y) \cdot \ind[\overline\CE] \big] &~\le~ p^2 \cdot \exp(2|\delta_1|) \cdot \E\big[\beta(x,y) \cdot \ind[\overline\CE] \big] \notag  \\
    &~\le~ \exp(2|\delta_1|) \cdot \E\big[\beta(x,y) \big],
\end{align}
since $\beta(x,y)$ is a non-negative random variable. Recall from \Cref{claim:second_moment_small_epsilon} that $\beta(x,y) \le \exp(\delta_1) \cdot \exp(Z)$ where $\delta_1 \le 1/\poly(d)$ and
\[ Z = d \epsilon^2 + 2\delta^2 \epsilon^2 r^2 + 2\delta^2 |\epsilon \langle Mx, My \rangle|.\]
Renormalizing $\overline{\epsilon} = \langle x,  n^{-1/2} y\rangle$ and $\overline{\theta} = \langle Mx, r^{-1} My \rangle$ and using that $\delta = \frac{\sqrt{d}}{\sigma \sqrt{n}}$ and $r \le  \sqrt{d}$, we have 
\begin{align*}
& Z ~\leq~ \left(\frac{d}{n} + \frac{2 d^2 }{\sigma^2 n^2} \right) \cdot \overline{\epsilon}^2 +  \frac{2 d \sqrt{d }}{\sigma^2 n \sqrt{n}} \cdot |\overline{\epsilon} \cdot \overline{\theta} | ~\leq~ (|\overline{\epsilon}| + |\overline{\theta}|)^2 / \lambda_{\min}^2 ,
\end{align*}
where we denote \[\lambda_{\min} = \frac{1}{3} \sqrt{\min\left\{ \frac{n}{d}, \frac{\sigma^2 n^2}{2d^2}, \frac{\sigma^2 n^{1.5}}{2d^{1.5}}\right\}}.\]
Note that $\lambda_{\min} = \omega_d(1) \cdot \sqrt{\log d}$ when $n =  \omega(d\log d)\sigma^{-4/3}$. 




We now bound the tails of $\overline\epsilon$ and $\overline{\theta}$, which will allow us to bound $\E[\exp(Z)]$.
Conditioned on any outcome of $y \sim \CG$, and as $\|y\| n^{-1/2}=1$, the second property  in
\Cref{lem:gram-schmidt} gives that   $\Pr_{x \sim \mathcal{G}} [ |\langle x,n^{-1/2}y \rangle | \geq t] \leq 2d^{C} \exp(-t^2/8)$.  Averaging over $y$ thus gives that
  \[ \Pr\left[|\overline\epsilon| \ge t\right]  \le 2 d^C \cdot e^{-t^2/8}. \]
 Similarly, as $\|My\| r^{-1} \leq 1$
 for any $y$ in the support of $\mathcal{G}$, we have that   $\Pr_{x \sim \mathcal{G}} [ |\langle Mx, r^{-1} My \rangle | \geq t] \leq 2d^{C} \exp(-t^2/8)$, and averaging over $y$ gives that
  \[
        \Pr\left[|\overline{\theta}| \ge t\right]
        \le 2 d^C\cdot e^{-t^2/8}. \] By a union bound, it follows that the random variable $X := |\overline\epsilon| + |\overline{\theta}|$ satisfies the tail condition of \Cref{lem:truncated_tail} with constant $C_1 = 2C$. So when $ \lambda_{\min}= \omega_d(1) \cdot \sqrt{\log d}$, the parameter $c_2 := \lambda_{\min}/ \sqrt{\log d}$ in \Cref{lem:truncated_tail} satisfies $32 C_1/c_2^2 = o_d(1)$. Therefore, \Cref{lem:truncated_tail} implies that
\begin{align*}
    \ \E[\beta(x,y)] &\le \exp(\delta_1) \cdot \E[\exp(Z)] \le (1+o_d(1))(1+o_d(1)) = 1+o_d(1). 
\end{align*}
Plugging the above in \eqref{eqn:1}, it follows that the second term 
\begin{align*}
    \E_{x,y \sim \CG}\big[P_x P_y \beta(x,y) \cdot \ind[\overline\CE] \big] &\le p^2 (1+o_d(1)).
\end{align*}
Combining this with \eqref{eqn:main} and \eqref{eqn:first}, we get that $\E_R[S^2] \le p^2 (1 + o_d(1))$. \qedhere
\end{proof}

We now prove the lemmas and claims used in the proof of \Cref{thm:main} above.

\subsection{Truncated Gram-Schmidt Distribution and Exponential Moments} \label{sec:Truncated}

\begin{proof}[Proof of \Cref{lem:gram-schmidt}]   

Consider running the Gram-Schmidt walk algorithm on the matrix $M$ stacked with the identity matrix, i.e. $\begin{pmatrix}M\\I_n\end{pmatrix}$, and let $\mathcal{G}_0$ be the distribution over colorings obtained as an output of the algorithm.

Since each column of the stacked matrix has Euclidean norm at most $2$, the properties of the Gram-Schmidt walk (\Cref{thm:hssz}) guarantees that $(x, Mx) \in \mathbb{R}^{n+d}$ where $x \in \{\pm 1\}^n$ and $Mx \in \mathbb{R}^d$ is 2-subgaussian. It follows that both $x$ and $Mx$ are $2$-subgaussian as well when $(x, Mx) \sim \CG_0$. 

To obtain a distribution $\CG$ where $\|M x\|_2$ is almost constant for each coloring $x \in \supp(\CG)$, we will truncate the distribution $\CG_0$ in such a way that the tails are also preserved up to $\poly(d)$ factors. Towards this end, we first note that with probability $1-e^{-cd}$, we have that $\|Mx\|_2 \le c'\sqrt{d}$ for constants $c$ and $c'$. This is because for any $\sigma$-subgaussian mean-zero random vector $X$, the Euclidean norm of $\|X\|$ has a subgaussian tail (e.g. 
Exercise 6.3.5 in \cite{V18}). In particular, $\p[\|X\|_2  \ge c_1 \sigma \sqrt{d} +  t] \le e^{-c_2 t^2/\sigma^2}$ for some universal constants $c_1, c_2 > 0$. Now, by a pigeonhole argument, for a large enough constant $C'$ there exist an annulus $W$ with width $\Delta = d^{-C'}$ and inner radius $r \le c'\sqrt{d}$ such that $\p_{x \sim \mathcal{G}_0}(x\in W) \ge d^{-C}$ for a constant $C$ depending on $C'$. 


We take the distribution $\mathcal{G}$ to be the probability measure of $\mathcal{G}_0$ conditioned on the event that $Mx \in W$. It then follows that for any coloring $x \sim \mathcal{G}$, we have $|\|Mx\|_2 - r| \leq \Delta$. Moreover, since $x$ and $Mx$ were $2$-subgaussian prior to conditioning, and the probability mass of the annulus is at least $d^{-C}$, conditioning can only increase the probability of any event by a factor of $d^{C}$. Thus, the tail bounds as stated in the statement of the lemma also follow. \qedhere
\end{proof}

\begin{proof}[Proof of \Cref{lem:truncated_tail}]
The assumption on $X$ implies that for any $t \geq 4 \cdot \sqrt{C_1 \log d}$, we have
\begin{align} \label{eq:tail_large}
    \p(X \geq t) \leq \exp(-t^2/16) . 
\end{align}
We express the expectation as an integration
\begin{align*} 
\E[\exp(X^2/\lambda^2)] ~ & = \int_0^\infty \p[\exp(X^2/\lambda^2) > s] ds  
\  = \int_0^\infty \p(X \geq \lambda \sqrt{\log s}) ds\\
& \leq 1 + c_3 + \int_{1 + c_3}^\infty \p(X \geq \lambda \sqrt{\log s}) ds.
\end{align*}
Let us set $c_3 = {32C_1}/c_2^2$, so that $c_3 \leq $ as $c_2 \geq \sqrt{32C_1}$.
For $s \geq 1 + c_3$, we have 
\[\lambda \sqrt{\log s} \geq c_2\sqrt{\log d} \cdot \sqrt{c_3/2} \ge 4 \cdot \sqrt{C_1 \log d}\] 
using that $\sqrt{\log (1+x)} \ge \sqrt{x/2}$ for $x \in [0,1]$ and as $c_3 \leq 1$.
So the condition $t \geq 4 \cdot \sqrt{C_1 \log d}$ for \eqref{eq:tail_large} is satisfied whenever $s \geq 1 + c_3$, and applying \eqref{eq:tail_large} to the above integration gives 
\begin{align*}
\int_{1 + c_3}^\infty \p(X \geq \lambda \sqrt{\log s}) ds &\leq \int_{1 + c_3}^\infty \exp(- \lambda^2 \log s/16) ds \\
& = \left(\frac{\lambda^2}{16}-1\right)^{-1} \cdot (1 + c_3)^{-\lambda^2/16 + 1} ~\leq~ \exp(-c_2^2 c_3 \log d/16). 
\end{align*}
By our choice of $c_3$, the above is at most $d^{-2C_1}$. Thus, it follows that $\E[\exp(X^2/\lambda^2)] \leq 1 + 32C_1/c_2^2 + d^{-2C_1}$. This proves the lemma. 
\end{proof}

\subsection{Proof of Claims from \Cref{sec:proof}} \label{sec:claims}
\label{sec:claim}

\begin{proof}[Proof of \Cref{claim:second_moment_small_epsilon}]
Since the rows of $R$ are independent, to compute the above ratio, it suffices to compute the ratio for a single row of $M + R$. Fix $i \in [d]$, and let $m = m_i$ and $r = r_i$ denote the $i^{\text{th}}$ row and define $a = a_i(x) := -m^\top x$ and $b = b_i(y):= -m^\top y$. We want to compare the ratio of $\p(r^\top x \in [a \pm \Delta], r^\top y \in [b \pm \Delta])$ to $\p(r^\top x \in [a \pm \Delta]) \cdot \p(r^\top y \in [b \pm \Delta])$.


Notice that $r^\top x$ and $r^\top y$ are Gaussian random variables with mean $0$, variance $1/\delta^2$, and covariance $\E_r[r^\top x r^\top y] ~=~ \E_r[x^\top r r^\top y] ~=~ \epsilon/\delta^2.$ Denoting the square $K :=  [a \pm \Delta] \times [b \pm \Delta]$, we have that $$\p(r^\top x \in [a \pm \Delta], r^\top y \in [b \pm \Delta]) = \mu_\epsilon(\delta K),$$
where $\mu_\epsilon$ is the 2-dimensional centered Gaussian measure with covariance matrix $\left(\begin{matrix}
1 & \epsilon \\
\epsilon & 1
\end{matrix}\right)$ and $\delta K$ denotes the $\delta$ scaling of $K$.  Similarly, we can write $\p(r^\top x \in [a \pm \Delta]) \cdot \p(r^\top y \in [b \pm \Delta]) = \mu(\delta K)$, where $\mu$ is the \emph{standard} 2-dimensional Gaussian measure.

We will compare the ratio $\mu_{\epsilon}(\delta K))/\mu(\delta K)$ by approximating the Gaussian measure over $\delta K$ with the density at the center and show the following bound
\begin{align}
\label{eqn:density}
    \ \mu_{\epsilon}(\delta K))/\mu(\delta K) \leq \exp \left(3\alpha + \epsilon^2 +   {\delta^2 \epsilon^2 (a^2 + b^2)} + 2{\delta^2 \epsilon a b}\right).
\end{align}

Since $(a_1(x), \ldots, a_d(x)) = Mx$ and $(b_1(x), \ldots, b_d(x)) = My$, using the above bound for all the rows $i \in [d]$, we have
\begin{align*} 
\frac{P_{x,y}}{P_x P_y} \leq \exp\left(3d\alpha + d \epsilon^2 + {\delta^2 \epsilon^2} \cdot (\|Mx\|_2^2 + \|My\|_2^2) + {\delta^2 \epsilon} \cdot \langle Mx, My \rangle \right). 
\end{align*}
Since $\|Mx\|_2^2 \le d + \poly(1/d)$ for every $x \in \supp(\CG)$, taking $\delta_1 = 4d\alpha$ gives the statement of the claim. To finish the proof we prove \eqref{eqn:density} now.

Abusing notation and denoting by $\mu(s,t)$ and $\mu_\epsilon(s,t)$ the corresponding densities at $(s,t) \in \mathbb{R}^2$, we have the following explicit formula for the density $\mu_\epsilon$: 
\begin{align*}
\mu_\epsilon(s,t)  = \frac{1}{2\pi \sqrt{1-\epsilon^2}} \cdot  \exp\left(-\frac{s^2 +t^2 - 2\epsilon st}{2(1-\epsilon^2)} \right).
\end{align*}

Since the edge length of the square $\delta K$ is $2\delta \Delta$, whenever $|\epsilon| \leq 1/2$, a direct calculation with the densities shows that 
\begin{align*}
    \frac{\sup_{(s,t) \in \delta K} \mu(s,t)}{\inf_{(s, t) \in \delta K} \mu(s,t)} ~=~ \exp\big(2\delta^2\Delta(|a|+|b|)\big) ~\leq~ \exp\big(2 \delta^2 \Delta (|a|+ |b| +  2 \Delta)\big) ~\le~ \exp(\delta_1),
\end{align*}
and that
\begin{align*}
    \frac{\sup_{(s,t) \in \delta K} \mu_\epsilon(s,t)}{\inf_{(s,t) \in \delta K} \mu_\epsilon(s,t)} \leq \exp(4 \delta^2 \Delta (|a|+ |b| + 2 \Delta)) \le  \exp(2\delta_1),
\end{align*}
where $\delta_1$ is as defined in the claim.
It follows that whenever $|\epsilon| \leq 1/2$, we can use the density at the center of $K$ to obtain
\begin{align*}
\frac{\mu_\epsilon(K)}{\mu(K)} &\leq \exp(3\delta_1) \cdot  \frac{\mu_\epsilon(\delta a, \delta b)}{\mu(\delta a, \delta b)} 
= \frac{1}{\sqrt{1- \epsilon^2}} \cdot \exp \left(3\delta_1 +  \frac{\delta^2 \epsilon^2 (a^2 + b^2)}{2(1 - \epsilon^2)} + \frac{\delta^2 \epsilon a b}{1 - \epsilon^2}\right) \\
& \leq \exp \big(3\delta_1 + \epsilon^2 +   {\delta^2 \epsilon^2 (a^2 + b^2)} + 2{\delta^2 \epsilon a b}\big),
\end{align*}
thus proving \eqref{eqn:density}.
\end{proof}

\begin{proof}[Proof of \Cref{claim:large_epsilon}]
We have $P_x = \prod_{i \in [d]} \p[r_i^\top x \in [a_i \pm \Delta]]$. For any fixed $i \in [d]$, $r_i^\top x$ is distributed as $\mathcal{N}(0,1/\delta^2)$, so after scaling the quantity $\p[r_i^\top x \in [a_i \pm \Delta]] = \mu(\delta \cdot I)$ where $I=[a_i \pm \Delta]$ and $\mu$ is the standard Gaussian measure in $\mathbb{R}$. Analogous to the proof of \Cref{claim:second_moment_small_epsilon}, one can approximate the Gaussian density at any  at the point in $I$ by the center point $a$, and compute similarly to the proof of \Cref{claim:second_moment_small_epsilon} that 
\begin{align*}
P_x = \prod_{i \in [d]} \p[r_i^\top x \in [a_i \pm \Delta]] 
= \left(\frac{\delta \Delta}{\sqrt{2 \pi}}\right)^d \exp\left(\alpha_x - \frac{\delta^2 \|Mx\|_2^2}{2} \right),
\end{align*}
for some small error $|\alpha_x| \leq 2\delta^2 \Delta (\|Mx\|_1 + d \Delta)$. As $\|Mx\|_2  \in [r \pm \Delta]$ and $r=O(\sqrt{d})$, we have that $\|Mx\|_1 =O(d)$ and the statement of the claim follows for some $\delta_x \le |\alpha_x| + 1/\poly(d) \le 1/\poly(d)$.
\end{proof}

\section{Conclusion}

For the Koml\'os problem, as studied in this paper, Gaussian noise is a natural way to model a smoothed analysis setting since the input vectors have Euclidean norm at most one. One can wonder whether similar results can be obtained with more general noise models, for instance, Bernoulli or other discrete noise models. Such noise models are also more conducive for smoothed analysis in other discrepancy settings, such as for the Beck-Fiala problem. The weighted second moment approach used here can also handle Bernoulli noise when the number of vectors $n \gg d^2$ but the second moment becomes difficult to control when $n$ is smaller. It remains an interesting open problem to see if Bernoulli or other discrete noise models can be handled for the regime $n \gg d \log d$.

\begin{small}
\bibliographystyle{alpha}
\bibliography{fullbib.bib}
\end{small}

\end{document}